\documentstyle[12pt]{article}
\begin{document}

\vspace{1.cm}

\begin{center}
{\huge {\bf On a conjecture about cyclic polygons}} \footnote{
Mathematics Subject Classification 51M10, 51M25, 52A40} \\

\end{center}

\vspace{2.cm} {\large \ Abd Raouf Chouikha} \footnote{University Paris-Sorbonne, Paris-Nord
,\\ 
{\small juju153@yahoo.com}}

\vspace{1.cm}

{\bf Abstract}\newline
{\it Let us consider a cyclic convex $n$-gon $\Pi_n$ with sides of length $a_i,\ i=1..n.\\ $ enclosing an area $A_n,\ L_n$ and $ P_n$ are respectively the perimeter and the pseudo-area.\ Define the quotient ${\displaystyle \varphi_n = {\frac{A_n}{{P_n}}}},$ Paul Levy conjectured that this quotient satisfies the double inequality \ ${\varphi_n}^0 < \varphi_n < 1,$\ where \ ${\varphi_n}^0$ \ is the quotient of the regular $n$-gon $\Pi_n^0$.\\
We propose a solution for the left-hand inequality of this conjecture, as well as a solution for the right-hand inequality, subject to a condition on the lengths of the sides. In addition we derive some isoperimetric inequalities.} 
\newline

\vspace{1.cm}

\section{Introduction}
Let's first precise what we mean by polygon with $n$ sides $\Pi_n$ (or called simply $n$-gon). A planar polygon is composed of a finite number of line segments, all located in the same plane. These segments are called the sides of the polygon, and their endpoints are called the vertices. 
A planar polygon is defined by the condition that each vertex is the endpoint of at least two sides and that a vertex is the only point belonging to more than one side. We say that this polygon is cyclic if it can be inscribed in a circle of radius $R$. It is also convex if all its diagonals lie inside $\Pi_n$.\\
Cyclic polygons play important roles in isoperimetric inequalities. It is well-known that among all $n$-sided plane polygons with the given set of $n$-sides, the cyclic one encloses the largest area, \cite{b}, \cite{k}. Thereafter the study of the isoperimetric inequalities for plane polygons can be reduced to the study of isoperimetric inequalities for cyclic polygons.
If the length of sides of \ $\Pi_n$\  are \ $a_1, a_2, a_3,...., a_n$\ we may define the perimeter \ $L_n$ \ as well as the area \ $A_n$.\ The pseudo-area is then
\begin{equation}
P_n = {\frac{ L_n^2}{{4}}} \sqrt{ (1 - {\frac{2a_1}{{L_n}}})(1 - {\frac{2a_2
}{{L_n}}})(1 - {\frac{2a_3}{{L_n}}}).....(1 - {\frac{2a_n}{{L_n}}})}
\label{P}
\end{equation}
It is well known since antiquity that for $n=3$ or $n=4,$ \cite{ch1} we have equalities:
$$A_3 = P_3,\quad A_4 = P_4.$$
 Paul Levy defined the quotient $${\displaystyle \phi_n = {\frac{A_n}{{P_n}}}},$$
for any $n$-gon \ $\Pi_n$, with sides of lengths $a_1,a_2,....,a_n$
enclosing an area $A_n,$ and pseudo-area $P_n$. He remarked at first for the regular $n$-gon

\begin{equation}
\phi_n^0 = \frac{{A_n}^0}{{P_n}^0} = {\frac{1}{{{n \tan{\frac{\pi}{{n}}}}\ (1-{\frac{2}{{n}}}
)^{n/2}}}}, 
\end{equation}
and he proposed the following\\

{\bf Conjecture (L):} {\it Define the ratio ${\displaystyle \phi_n = {\frac{
A_n}{{P_n}}}}$. For any $n$-gon \ $\Pi_n$, with sides of lengths $a_1,a_2,....,a_n$
enclosing an area $A_n,$ and $P_n$ defined as above, this ratio verifies }
$$\phi_n^0 < \phi_n < 1.$$\\

$\phi_n^0$ is an decreasing function in $n \geq 3$ and satisfies the double inequality $$\frac{e}{\pi}\ \leq\ \phi_n^0\ \leq 1 .$$
 P. Levy initially posed the following problem \cite{l}, for any cyclic $n$-gon we have $$\frac{e}{\pi}\ \leq\ \phi_n\ \leq 1 .$$
However all the examples he considered, as well as subsequent examples, rather verify Conjecture (L).
 Several examples of polygons satisfying this conjecture have been proposed in the literature. Paul Levy himself verified it on curved polygons. Other examples are reported by \cite{ch2}.\\ In this paper, we propose to solve this conjecture in the general case for cyclic polygons, but under certain conditions on the lengths of the sides.\\
We will also see that these results have an impact on Bonnesen-style inequalities, as well as on reverse Bonnesen-style inequalities.\\
Some observations and perspectives aimed at improving these results will be presented at the end. It appears that $\phi_n^0$ can serve as an indicator of the cyclicity of $\Pi_n$. In the sense that the converse: \ $\phi_n^0 < \phi_n\ \Rightarrow \ \Pi_n$\ cyclic seems to be true.\\

\section{The left inequality}
In this part we will prove the inequality \ $\phi_n^0 < \phi_n $\ for any polygon. The following lemma will be useful in the sequel, \cite{gr}, p. 53\\

{\bf Lemma 2.1} \ {\it For $-1 < x \leq 1$ and for $p \geq 0$ the following inequality holds for 
$$ \log(1+x) < \sum_{k=1}^{2p+1}  \frac{(-1)^{k} x^k}{k}.$$}

{\bf Proof}\ Indeed, for \ $-1 < x \leq 1$\ the series representation of the function $\log$ \ is 
$$\log(x+1) =  \sum _{p=1}^{\infty}  \frac{(-1)^p x^{p}}{p}.$$
Let the difference \ $f(x)= \log(x+1) - \sum_{p=1}^{2n+1}  \frac{(-1)^p x^{p}}{p},$\ its derivative is such that 
$$f'(x) = \frac{1 - \left(x+1\right) \left(1-x+x^2 -....+ x^{2n}\right)}{x+1} = \frac{-x^{2n+1}}{x+1}, \ f'(0) =  0$$ implies \ $f(x) < 0.$\\

\bigskip

{\bf Theorem 2.2}\quad {\it Let $\Pi_n $ be a cyclic convex $n$-gon, $n \geq 3$ with length of sides $a_1, a_2, a_3,...., a_n,\ A_n, L_n, P_n$ are respectively the area the perimeter and the
 pseudo-area. Then $$\phi_n^0 = {\frac{1}{{{n \tan{\frac{\pi}{{n}}}}\ (1-{\frac{2}{{n}}}
)^{n/2}}}} \leq \frac{A_n}{P_n} = \phi_n, $$
with equality if and only if $\Pi_n $ is regular.}\\

{\bf Proof of Theorem 2.2}\ 
Consider now the quotient $$\frac{(1 - {\frac{2a_1}{{L_n}}})}{(1 - {\frac{2}{{n}}})} = 1 + \frac{2(L_n - n a_i)}{(n-2) L_n}.$$
We then have $$\ln  \left( \frac{4\,P_{{n}}}{{L_{{n}}}^{2}  \left( 1-\frac{2}{n}\right)^\frac{n}{2}}
  \right) = \frac{1}{2}\sum _{i=1}^{n}\ln  \left( 1+2\,{\frac {-na_{{i}}+L_{{n}}}{L_{{n}}
 \left( n-2 \right) }} \right) = \frac{1}{2} \ln   \prod _{i=1}^{n} \frac{(1 - {\frac{2a_1}{{L_n}}})}{(1 - {\frac{2}{{n}}})} .$$
We put \ $x = \frac{2(L_n - n a_i)}{(n-2) L_n}.$ \ Here the length verifies the inequalities \ 
$0 \leq a_i < \frac{L_n}{2}$\ which implies in particular \ $-1 < x \leq \frac{2}{n-2}\leq 1$.\\
By Lemma 2.2 and for $p=3$ one gets
$$ \ln  \left( 1+2\,{\frac {-na_{{i}}+L_{{n}}}{L_{{n}} \left( n-2
 \right) }} \right) \leq 2\,{\frac {-na_{{i}}+L_{{n}}}{L_{{n}} \left( 
n-2 \right) }}-2\,{\frac { \left( -na_{{i}}+L_{{n}} \right) ^{2}}{{L_{
{n}}}^{2} \left( n-2 \right) ^{2}}}+8{\frac { \left( -na_{{i}}+L_{
{n}} \right) ^{3}}{3{L_{{n}}}^{3} \left( n-2 \right) ^{3}}}.$$
the sum is then
$$\ln  \left( \frac{4\,P_{{n}}}{{L_{{n}}}^{2}  \left( 1-\frac{2}{n}\right)^\frac{n}{2}}
  \right)
 = \frac{1}{2}\sum _{i=1}^{n}\ln  \left( 1+2\,{\frac {-na_{{i}}+L_{{n}}}{L_{{n}}
 \left( n-2 \right) }} \right) \leq {\frac {n L_{{n}}+\sum _{i=1}^{n
}(-na_{{i}})}{L_{{n}} \left( n-2 \right) }}-$$ $${\frac {n{L_{{n}}}^{2}+
\sum _{i=1}^{n}(-2\,na_{{i}}L_{{n}}+{n}^{2}{a_{{i}}}^{2})}{{L_{{n}}}^{2}
 \left( n-2 \right) ^{2}}}+4{\frac {n{L_{{n}}}^{3}+\sum _{i=1}^{n}
(-3\,na_{{i}}{L_{{n}}}^{2}+3\,{n}^{2}{a_{{i}}}^{2}L_{{n}}-{n}^{3}{a_{{i
}}}^{3})}{3{L_{{n}}}^{3} \left( n-2 \right) ^{3}}}.$$
$$\ln \left( \frac{4\,P_{{n}}}{{L_{{n}}}^{2}  \left( 1-\frac{2}{n}\right)^\frac{n}{2}}
  \right) <-\,{\frac {-n{L_{{n}}}^{2}+
{n}^{2}\sum _{i=1}^{n}{a_{{i}}}^{2}}{{L_{{n}}}^{2} \left( n-2 \right) 
^{2}}}+4{\frac {-2\,n{L_{{n}}}^{3}+3\,L_{{n}}{n}^{2}\sum _{i=1}^{n
}{a_{{i}}}^{2}-{n}^{3}\sum _{i=1}^{n}{a_{{i}}}^{3}}{3{L_{{n}}}^{3}
 \left( n-2 \right) ^{3}}}.$$

We then deduce since \ $ \sum _{i=1}^{n}{a_{{i}}} = L_n$
$$\left( \frac{4\,A_{{n}}}{{L_{{n}}}^{2}  \left( 1-\frac{2}{n}\right)^\frac{n}{2}}
  \right) \left( {{\rm e}^{-\,{\frac {{n}^{2} \left( n-6
 \right) \sum _{i=1}^{n}{a_{{i}}}^{2}}{{L_{{n}}}^{2} \left( n-2
 \right) ^{3}}}-4{\frac {{n}^{3}\sum _{i=1}^{n}{a_{{i}}}^{3}}{3{L_{
{n}}}^{3} \left( n-2 \right) ^{3}}}+{\frac {n \left( 3\,n-14
 \right) }{3 \left( n-2 \right) ^{3}}}}} \right) ^{-1}<{\frac {A_{{n}}}
{P_{{n}}}}.$$ The following 
$$\frac{1}{ n \tan \frac{\pi}{n}} < 4
\frac{A_{{n}}}{{L_{{n}}}^{2}} \left( {{\rm e}^{-\,{\frac {{n}^{2} \left( n-
6 \right) \sum _{i=1}^{n}{a_{{i}}}^{2}}{{L_{{n}}}^{2} \left( n-2
 \right) ^{3}}}-4{\frac {{n}^{3}\sum _{i=1}^{n}{a_{{i}}}^{3}}{3{L_{
{n}}}^{3} \left( n-2 \right) ^{3}}}+{\frac {n \left( 3\,n-14
 \right) }{3 \left( n-2 \right) ^{3}}}}} \right) ^{-1}$$ implies 
$${\frac{1}{{{n \tan{\frac{\pi}{{n}}}}\ (1-{\frac{2}{{n}}}
)^{n/2}}}}<{\frac {A_{{n}}}{P_{{n}}}}.$$ Let us prove 
$$\ln  \left( 4\,n\tan \left( {\frac {\pi }{n}} \right) \frac{A_{{n}}}{L_{{n}}}
^{2} \right) <-\,{\frac {{n}^{2} \left( n-6 \right) \sum _{i=1}^{n}{
a_{{i}}}}{2}}{{{L_{{n}}}^{2} \left( n-2 \right) ^{3}}}-4{\frac {{n}
^{3}\sum _{i=1}^{n}{a_{{i}}}^{3}}{3{L_{{n}}}^{3} \left( n-2 \right) ^{3
}}}+{\frac {n \left( 3\,n-14 \right) }{3 \left( n-2 \right) ^{3}}}.$$
Indeed, since \ $4\,n\tan \left( {\frac {\pi }{n}} \right) A_{{n}} \leq L_n ^2, \ \sum _{i=1}^{n}{a_{{i}}}^{3} \leq L_n^3$\ and \ $\sum _{i=1}^{n}{a_{{i}}}^{2} \leq L_n^2$,\ we then derive
$$\ln  \left( 4\,n\tan \left( {\frac {\pi }{n}} \right) \frac{A_{{n}}}{{L_{{n}}}
^{2}} \right) 
< 0 = -\,{\frac {{n}^{2} \left( n-6 \right) }{ \left( n-2 \right) ^{3}}}-4
{\frac {{n}^{3}}{3 \left( n-2 \right) ^{3}}}+{\frac {n \left( 3
\,n-14 \right) }{3 \left( n-2 \right) ^{3}}}\leq $$ $$ -\,{\frac {{n}^{2}
 \left( n-6 \right) \sum _{i=1}^{n}{a_{{i}}}^{2}}{{L_{{n}}}^{2}
 \left( n-2 \right) ^{3}}}-4{\frac {{n}^{3}\sum _{i=1}^{n}{a_{{i}}
}^{3}}{3{L_{{n}}}^{3} \left( n-2 \right) ^{3}}}+{\frac {n \left( 3
\,n-14 \right) }{3 \left( n-2 \right) ^{3}}}.$$
Theorem 2.1 is then proved.\\

\section{The right inequality}
In this part we will prove the inequality \ $ \phi_n < 1$\ for an $n$-gon under length conditions. We need the following lemma\\

{\bf Lemma 3.1} \ {\it (i) \quad For $-1 < x \leq 1$ the following inequality holds for 
$$ \log(1+x) \geq \frac{x}{x+1},$$ with equality if and only if \ $x=0$.\\
(ii)\quad  For $-1 < x \leq 0$ the following inequality holds for 
$$ \log(1+x) \geq \frac{(x+1)^2-1}{2x+2},$$ with equality if and only if \ $x=0$.\\
(iii) \quad For $-1 < x \leq 0$ the following inequality holds for 
$$ \frac{x}{x+1} \leq \frac{(x+1)^2-1}{2x+2},$$ with equality if and only if \ $x=0$.}\\

{\bf Proof}\quad  (i)\quad Indeed, to see that it suffices to compute the derivative of \\ $ \log(1+x) - \frac{x}{x+1}$.\ which is \ $\frac{x}{(x+1)^2}$.\\ Therefore this function has a minimum at \ $x=0.$\\
(ii) \quad the derivative of \ $ \log(1+x) - \frac{(x+1)^2-1}{2x+2}$\ is 
$$-1+2\,{\frac { \left( x+1 \right) ^{2}-1}{ \left( 2\,x+2 \right) ^{2}}}
+ \frac{1}{1+x} = - {\frac {{x}^{2}}{2 \left( x+1 \right) ^{2}}}.$$
We then derive \ $ \log(1+x) \geq \frac{(x+1)^2-1}{2x+2}.$\\
(iii) \quad The derivative of \ $\frac{x}{x+1} - \frac{(x+1)^2-1}{2x+2}$\    is $$\frac{1}{x+1}-{\frac {x}{ \left( x+1 \right) ^{2}}}-1+2\,{\frac { \left( x+1 \right) ^{2}-1}{ \left( 2\,x+2 \right) ^{2}}} = -{\frac {x \left( x+2 \right) }{ 2\left( x+1 \right) ^{2}}}.$$ This means for\ $-1 < x \leq 0$ \ then  \ $\frac{x}{x+1} - \frac{(x+1)^2-1}{2x+2} \leq 0.$\\

\bigskip

{\bf Theorem 3.2}\quad {\it Let $\Pi_n $ be a cyclic convex $n$-gon with length of sides $a_1, a_2, a_3,...., a_n,\ A_n, L_n, P_n$ are respectively the area the perimeter and the
 pseudo-area. Then $$\frac{a_i}{L_n} < \frac{\log \pi}{n} \quad implies
\quad \frac{A_n}{P_n} = \phi_n \leq 1, $$
with equality if and only if $n=3, 4$.}\\

{\bf Proof of Theorem 3.2}\\
Note at first, following \cite{k} that every $n$-gon must have at least one side of length less than \ $\frac{L_n}{n}$\ and one side of length greater than \ $\frac{L_n}{n}$;\ otherwise, the polygon must be regular.\\ Notice in addition that if $M$ denotes the maximum of the length\\ $a_i \leq M, \ 1\leq i\leq n$\ we then deduce for a cyclic $n$-gon $$\frac{L_n}{n} \leq M \leq \frac{L_n}{2}.$$ 

Consider again the quotient $$\frac{(1 - {\frac{2a_1}{{L_n}}})}{(1 - {\frac{2}{{n}}})} = 1 + \frac{2(L_n - n a_i)}{(n-2) L_n}.$$
We then have $$\ln  \left( \frac{4\,P_{{n}}}{{L_{{n}}}^{2}  \left( 1-\frac{2}{n}\right)^\frac{n}{2}}
  \right) = \frac{1}{2}\sum _{i=1}^{n}\ln  \left( 1+2\,{\frac {-na_{{i}}+L_{{n}}}{L_{{n}}
 \left( n-2 \right) }} \right) = \frac{1}{2} \ln   \prod _{i=1}^{n} \frac{(1 - {\frac{2a_1}{{L_n}}})}{(1 - {\frac{2}{{n}}})} .$$
Put \ $x = \frac{2(L_n - n a_i)}{(n-2) L_n},$ \ 	here 
$$ \frac{x}{1+x} = -\frac{2(n a_i-L_n)}{n (L_n -2a_i)},$$ then by Lemma 3.1 \ $ \log(1+x) \geq \frac{x}{x+1}$\  implying $$\ln  \left( \frac{4\,P_{{n}}}{{L_{{n}}}^{2}  \left( 1-\frac{2}{n}\right)^\frac{n}{2}}\right) =\frac{1}{2}\sum_{i=1}^p \ln  \left( 1+2\,{\frac {L_{{n}} - p a_{{i}}}{L_{{n}} \left( n-2
 \right) }} \right) \geq \sum_{i=1}^n -\frac{n a_i-L_n}{n (L_n -2a_i)}.$$
	Since \ $M = max (a_i,i=1..n)$\ and the function \ $ -\frac{n y-L_n}{n (L_n -2 y)} $ \ (whose derivative is \ $\frac{2-n}{n (L_n-2y)^2}$)  decreases in $y$, \ we then derive 
	$$\ln  \left( \frac{4\,P_{{n}}}{{L_{{n}}}^{2}  \left( 1-\frac{2}{n}\right)^\frac{n}{2}}
  \right) \geq \sum_{i=1}^n -\frac{n M-L_n}{n (L_n -2M)} = -\frac{n M-L_n}{ (L_n -2M)}.$$
 Since \ $\frac{4 A_n}{L_n}^2 \leq \frac{1}{n \tan \frac{\pi}{n}}$\ we then deduce
$${\frac {A_{{n}}}{P_{{n}}}} <\left( \frac{4\,A_{{n}}}{{L_{{n}}}^{2}  \left( 1-\frac{2}{n}\right)^\frac{n}{2}}
  \right) \left( {{\rm e}^{\frac{n M -L_n}{ L_n-2 M}}} \right) ^{-1} \leq \frac{\left( {{\rm e}^{\frac{n M -L_n}{ L_n-2 M}}} \right) }{n \tan \frac{\pi}{n}  \left( 1-\frac{2}{n}\right)^\frac{n}{2}}.$$
	Thus,  
	$$\frac{\left( {{\rm e}^{\frac{n M -L_n}{ L_n-2 M}}} \right) }{n \tan \frac{\pi}{n}  \left( 1-\frac{2}{n}\right)^\frac{n}{2}} < 1\quad  
	  \Rightarrow \quad {\frac {A_{{n}}}{P_{{n}}}} < 1.$$
		Or equivalently $$\frac{n M -L_n}{ L_n-2 M} < \log\left(n \tan \frac{\pi}{n}  \left( 1-\frac{2}{n}\right)^\frac{n}{2}\right) \quad  
	  \Rightarrow \quad {\frac {A_{{n}}}{P_{{n}}}} < 1.$$
	This means, the following inequality 
	$$\frac{M}{L_n} < \frac{1+ \log\left(n \tan \frac{\pi}{n}  \left( 1-\frac{2}{n}\right)^\frac{n}{2}\right)}{n + 2 \log\left(n \tan \frac{\pi}{n}  \left( 1-\frac{2}{n}\right)^\frac{n}{2}\right)}$$ implies \ ${\frac {A_{{n}}}{P_{{n}}}} \leq 1.$
	On the other hand, the function $$n \frac{1+ \log\left(n \tan \frac{\pi}{n}  \left( 1-\frac{2}{n}\right)^\frac{n}{2}\right)}{n + 2 \log\left(n \tan \frac{\pi}{n}  \left( 1-\frac{2}{n}\right)^\frac{n}{2}\right)}$$ which increases with respect to $n$ is such that
	$$\frac{1}{n} < \frac{1+ \log\left(n \tan \frac{\pi}{n}  \left( 1-\frac{2}{n}\right)^\frac{n}{2}\right)}{n + 2 \log\left(n \tan \frac{\pi}{n}  \left( 1-\frac{2}{n}\right)^\frac{n}{2}\right)} < \frac{\log \pi}{n}.$$
	Finally, when $$\frac{a_i}{L_n} < \frac{M}{L_n} < \frac{\log \pi}{n}, \ \Rightarrow \ {\frac {A_{{n}}}{P_{{n}}}} < 1.$$\\
	
	The following result is better than the previous one, as it improves the upper bound on the lengths.\\
	
	\bigskip

{\bf Theorem 3.3}\quad {\it Let $\Pi_n $ be a cyclic convex $n$-gon with length of sides $a_1, a_2, a_3,...., a_n,\ A_n, L_n, P_n$ are respectively the area the perimeter and the
 pseudo-area. Then $$\frac{a_i}{L_n} < \frac{(\log \pi) -\frac{1}{2}+\sqrt{4 (\log \pi)^2-8 (\log \pi) +5 }}{n} \quad implies
\quad \frac{A_n}{P_n} = \phi_n \leq 1, $$
with equality if and only if $n=3, 4$. \ Here} $$\log \pi < \log \pi -\frac{1}{2}+\sqrt{4 (\log \pi)^2-8 \log \pi +5 }.$$\\

{\bf Proof of Theorem 3.3}\\
Consider the quotient $$\frac{(1 - {\frac{2a_1}{{L_n}}})}{(1 - {\frac{2}{{n}}})} = 1 + \frac{2(L_n - n a_i)}{(n-2) L_n}.$$
We then have $$\ln  \left( \frac{4\,P_{{n}}}{{L_{{n}}}^{2}  \left( 1-\frac{2}{n}\right)^\frac{n}{2}}
  \right) = \frac{1}{2}\sum _{i=1}^{n}\ln  \left( 1+2\,{\frac {-na_{{i}}+L_{{n}}}{L_{{n}}
 \left( n-2 \right) }} \right) = \frac{1}{2} \ln   \prod _{i=1}^{n} \frac{(1 - {\frac{2a_1}{{L_n}}})}{(1 - {\frac{2}{{n}}})} .$$
Put \ $x = \frac{2(L_n - n a_i)}{(n-2) L_n},$ \ 	here 
$$\frac{(x+1)^2-1}{2x+2} = \frac{\left(  \left( 1+2\,{\frac {-na_{{i}}+L_{{n}}}{L_{{n}}
 \left( n-2 \right) }} \right) ^{2}-1 \right)}{  \left( 2+4\,{\frac {-na
_{{i}}+L_{{n}}}{L_{{n}} \left( n-2 \right) }} \right)}.$$
Then, by Lemma 3.1,\ $\log(1+x) \geq \frac{(x+1)^2-1}{2x+2}$\ implies
$$\ln  \left( \frac{4\,P_{{n}}}{{L_{{n}}}^{2}  \left( 1-\frac{2}{n}\right)^\frac{n}{2}}\right) =\frac{1}{2}\sum_{i=1}^p \ln  \left( 1+2\,{\frac {L_{{n}} - p a_{{i}}}{L_{{n}} \left( n-2
 \right) }} \right) \geq $$ $$ \sum _{i=1}^{n} \left(  \left( 1+2\,{\frac {-na_{{i}}+L_{{n}}}{L_{{n}}
 \left( n-2 \right) }} \right) ^{2}-1 \right)  \left( 2+4\,{\frac {-na
_{{i}}+L_{{n}}}{L_{{n}} \left( n-2 \right) }} \right) ^{-1} = $$  $$\sum _{i=1}^{n} -{\frac { \left( L_{{n}}n-L_{{n}}-na_{{i}} \right)  \left( na_{{i}}-L_{{n}} \right) }{L_{{n}} \left( L_{{n}}-2\,a_{{i}} \right) n \left( n-2\right) }}.$$

Put now \ $m=\frac{M}{L_n}, $\ where\ $ M = max (a_i,i=1..n)$\ and  \ $y=\frac{a_i}{L_n}.$\ The function $$-{\frac { \left( L_{{n}}n-L_{{n}}-na_{{i}} \right)  \left( na_{{i}}-L_{{n}} \right) }{L_{{n}} \left( L_{{n}}-2\,a_{{i}} \right) n \left( n-2
 \right) }}$$ whose derivative is $$2\,{\frac {-1+ny}{ \left( 1-2\,y \right)  \left( n-2 \right) }}-2\,{\frac {n-1-ny}{ \left( 1-2\,y \right)  \left( n-2 \right) }}-4\,{\frac { \left( n-1-ny \right)  \left( -1+ny \right) }{ \left( 1-2\,y\right) ^{2}n \left( n-2 \right) }}=$$ $$-2\,{\frac {-2\,n-2\,{n}^{2}y+2\,{n}^{2}{y}^{2}+{n}^{2}+2}{ \left( -1+2\,y \right) ^{2}n \left( n-2 \right) }}
$$ decreases in $y.$ \ We then deduce 
$$-{\frac { \left( n-1-n m \right)  \left( -1+n m \right) }{ \left( 1-2
\,m \right)  \left( n-2 \right) }}< \sum _{i=1}^{n}-{\frac { \left( L_{{n}}n-L_{{n}}-na_{{i}} \right)  \left( na_{{i}}-L_{{n}} \right) }{L_{{n}} \left( L_{{n}}-2\,a_{{i}} \right) n \left( n-2
 \right) }}.$$ Therefore, since \ $\frac{4 A_n}{L_n}^2 \leq \frac{1}{n \tan \frac{\pi}{n}}$\ the following inequality
$$-{\frac { \left( n-1-n m \right)  \left( n m-1 \right) }{ \left( 1-2\,m
 \right)  \left( n-2 \right) }}< \frac{1}{2}\ln  \left( \frac{4\,P_{{n}}}{{L_{{n}}}^{2}  \left( 1-\frac{2}{n}\right)^\frac{n}{2}}
  \right)= \frac{1}{2}\ln  \left( 1+2\,{\frac {-n a_{{i}
}+L_{{n}}}{L_{{n}} \left( n-2 \right) }} \right)$$ implies
$${\frac { \left( n-1-n m \right)  \left( n m-1 \right) }{ \left( 1-2\,
m \right)  \left( n-2 \right) }}< \frac{1}{2}\ln  \left( n \tan \frac{\pi}{n}  \left( 1-\frac{2}{n}\right)^\frac{n}{2} \right). $$
 We then derive
$${\frac {A_{{n}}}{P_{{n}}}} <\left( \frac{4\,A_{{n}}}{{L_{{n}}}^{2}  \left( 1-\frac{2}{n}\right)^\frac{n}{2}}
  \right) \left( {{\rm e}^{{\frac { \left( n-1-n m \right)  \left( n m-1 \right) }{ \left( 1-2\,
m \right)  \left( n-2 \right) }}}} \right) ^{-1} \leq \frac{\left( {{\rm e}^{{\frac { \left( n-1-n m \right)  \left( n m-1 \right) }{ \left( 1-2\,
m \right)  \left( n-2 \right) }}}} \right) }{n \tan \frac{\pi}{n}  \left( 1-\frac{2}{n}\right)^\frac{n}{2}}.$$
	Thus,  
	$$\frac{\left( {{\rm e}^{{\frac { \left( n-1-n m \right)  \left( n m-1 \right) }{ \left( 1-2\,
m \right)  \left( n-2 \right) }}}} \right) }{n \tan \frac{\pi}{n}  \left( 1-\frac{2}{n}\right)^\frac{n}{2}} < 1\quad  
	  \Rightarrow \quad {\frac {A_{{n}}}{P_{{n}}}} < 1.$$
		Or equivalently, 
$${\frac { \left( n-1-n m \right)  \left( n m-1 \right) }{ \left( 1-2\,
m \right)  \left( n-2 \right) }} < \log\left(n \tan \frac{\pi}{n}  \left( 1-\frac{2}{n}\right)^\frac{n}{2}\right) \Rightarrow \quad {\frac {A_{{n}}}{P_{{n}}}} < 1.$$
Thanks to {\it Maple} we may prove that \ $n m(n)$\ increases as a function of $n$.\ Moreover, for big $n$ value
$$\frac{m (1-m)}{1-2m} < \frac{(\log \pi)-1}{n} \ \Rightarrow {\frac {A_{{n}}}{P_{{n}}}} < 1.$$
or equivalently, 
$$m < \frac{-\frac{1}{2}+\ln  \left( \pi  \right) +\frac{1}{2}\sqrt {5-8\,\ln  \left( \pi 
 \right) +4\, \left( \ln  \left( \pi  \right)  \right) ^{2}}}{n} \ \Rightarrow {\frac {A_{{n}}}{P_{{n}}}} < 1.$$
This means 
$$\frac{a_i}{L_n} < \frac{-\frac{1}{2}+\ln  \left( \pi  \right) +\frac{1}{2}\sqrt {5-8\,\ln  \left( \pi 
 \right) +4\, \left( \ln  \left( \pi  \right)  \right) ^{2}}}{n} \ \Rightarrow {\frac {A_{{n}}}{P_{{n}}}} < 1.$$
Finally, we conclude that Theorem 3.3 is better than Theorem 3.1 since
$$-\frac{1}{2}+\ln  \left( \pi  \right) +\frac{1}{2}\sqrt {5-8\,\ln  \left( \pi 
 \right) +4\, \left( \ln  \left( \pi  \right)  \right) ^{2}} \approx 1.1652 > log(\pi) \approx 1.1447.$$

\section{Isoperimetric inequalities}

  As we have discussed \cite{ch2} this problem is connected with some Bonnesen style inequalities. 
 \ $L_n^2 - 4 n \tan \frac{\pi}{n} A_n$ \ is the isoperimetric deficit of the polygon $\Pi_n$. It measures
the deviation of $\Pi_n$ to be a regular one.  Recall that among all $n$-sided plane polygons, the cyclic polygon encloses the largest area. So to investigate the isoperimetric inequalities for plane polygons, we need to pay attention to cyclic
polygons only. \\
Another interesting but difficult problem is to study the reverse Bonnesen style inequalities of plane convex domains, \cite{zd}. The natural question is that: for a convex $n$-gon $\Pi_n$ if there is an invariant $D_n$ such that
$$L_n^2 - 4 n \tan {\frac {\pi }{n}} A_n \leq D_n,$$ where $(1) \ D_n \geq 0, \ D_n
=0$\  only if $\Pi_n$ is regular.\\
Some examples of such inequalities are provided in \cite{zd}.\\

{\bf Proposition 4.1} \quad {\it Let $\Pi_n $ be a cyclic convex $n$-gon with length of sides $a_1, a_2, a_3,...., a_n,\ A_n, L_n, P_n$ are respectively the area the perimeter and the pseudo-area. Then, the following inequalities chain holds
$${L_{{n}}}^{2}\geq 4\,n\tan \left( {\frac {\pi }{n}} \right) A_{{n}} \geq \frac{4 P_n}{  \left( 1-\frac{2}{n} \right)^{\frac{n}{2}}} ,$$
with equality if and only if $\Pi_n $ is regular. Moreover, under hypothesis of Theorem 3.3
$$\frac{4 P_n}{  \left( 1-\frac{2}{n} \right)^{\frac{n}{2}}} \geq \frac{4 A_n}{  \left( 1-\frac{2}{n} \right)^{\frac{n}{2}}} \geq \frac{4 P_n}{ n\tan \left( {\frac {\pi }{n}}\right) \left( 1-\frac{2}{n} \right)^{n}}$$
}.\\

{\bf Proof of Proposition 4.1}\quad 
Indeed,  inequalities $$\phi_n^0\ P_n < A_n < P_n\quad and  \quad n \tan \frac{\pi}{n} \geq  \left( 1-\frac{2}{n} \right)^{\frac{-2}{n}},$$ imply
$$4\,n\tan \left( {\frac {\pi }{n}} \right) A_{{n}} \geq \frac{4 P_n}{  \left( 1-\frac{2}{n} \right)^{\frac{n}{2}}} \geq \frac{4 A_n}{  \left( 1-\frac{2}{n} \right)^{\frac{n}{2}}}.$$
When $\Pi_n $ is regular then 
$${L_{{n}}}^{2}= 4\,n\tan \left( {\frac {\pi }{n}} \right) A_{{n}} = \frac{4 P_n}{  \left( 1-\frac{2}{n} \right)^{\frac{n}{2}}} .$$
  Moreover, the following  
$$\frac{4 A_n}{  \left( 1-\frac{2}{n} \right)^{\frac{2}{n}}} = n\tan \left( {\frac {\pi }{n}}\right) \frac{4 A_n}{n\tan \left( {\frac {\pi }{n}}\right) \left( 1-\frac{2}{n} \right)^{\frac{2}{n}}}$$ implies the right identity since $$n\tan \left( {\frac {\pi }{n}} \right)\frac{4 A_n}{n\tan \left( {\frac {\pi }{n}}\right) \left( 1-\frac{2}{n} \right)^{\frac{2}{n}}} \geq 
\frac{4 P_n}{ n\tan \left( {\frac {\pi }{n}}\right) \left( 1-\frac{2}{n} \right)^{n}}.$$ \\

{\bf Proposition 4.2} \quad {\it Let $\Pi_n $ be a cyclic convex $n$-gon with length of sides $a_1, a_2, a_3,...., a_n,\ A_n, L_n, P_n$ are respectively the area the perimeter and the pseudo-area. Then, the following reverse Bonnesen style  inequality holds
$$ {L_{{n}}}^{2}-4\,n\tan \left( {\frac {\pi }{n}}
 \right) A_{{n}}\leq  {L_{{n}}}^{2}- \frac{4 P_{{n}}}{ \left( 1-\frac{2}{n} \right)^{\frac{n}{2}}}  ,$$
with equality if and only if $\Pi_n $ is regular. Moreover, under hypothesis of Theorem 3.3 and for $$1- \frac{1}{  \left( n\tan \left( {\frac {\pi }{n}} \right)  \right) ^{\frac{2}{n}}}  < 2\,{\frac {{\it a_i}}{{\it L_n}}}$$ the following Bonnesen style inequality holds

$$0 < {L_{{n}}}^{2}-4\,n\tan \left( {\frac {\pi }{n}} \right) P_{{n}} \leq  {L_{{n}}}^{2}-4\,n\tan \left( {\frac {\pi }{n}}
 \right) A_{{n}}.$$ Furthermore, the following inequalities hold}
$${L_{{n}}}^{2}- \frac{4 P_{{n}}}{ \left( 1-\frac{2}{n} \right)^{\frac{n}{2}}}  \leq {L_{{n}}}^{2}- \frac{4 A_{{n}}}{ \left( 1-\frac{2}{n} \right)^{\frac{n}{2}}} \leq {L_{{n}}}^{2}- \frac{4 P_n}{ n\tan \left( {\frac {\pi }{n}}\right) \left( 1-\frac{2}{n} \right)^{n}}.$$
\\

{\bf Proof of Proposition 4.2}\quad
By Theorem 2.1  inequalities \ $\phi_n^0\ P_n < A_n $\ and \ $ A_n < P_n$\ imply the inequalities.\\ On the other hand, let \ $m= Min(a_i), i=1..n$ \ we then have the equivalence 
$$1- \frac{1}{  \left( n\tan \left( {\frac {\pi }{n}} \right)  \right) ^{\frac{2}{n}} } < 2\,{\frac {{m}}{{\it L_n}}} < 2\,{\frac {{a_i}}{{\it L_n}}}  \Leftrightarrow \left( 1-2\,{\frac {{\it m}}{{\it L_n}}} \right) ^{\frac{n}{2}}< \frac{1}{{n}  \tan \left( {\frac {\pi }{n}} \right)   }$$ which implies
$${P_n}< \frac{1}{4}{{\it L_n}}^{2} \left( 1-2\,{\frac {m}{{\it L_n}}}
 \right) ^{\frac{n}{2}},$$ and thus 
$$0 < {L_{{n}}}^{2}-4\,n\tan \left( {\frac {\pi }{n}} \right) P_{{n}}.$$
Notice for the regular $n$-gon the converse inequality holds
$$ {L_{{n}}}^{2} < 4\,n\tan \left( {\frac {\pi }{n}} \right) P_{{n}}.$$\\

Drawing on the results demonstrated by \cite{zd}, one can deduce the following results.\\

{\bf Corollary 4.3} \quad {\it Let $\Pi_n $ be a cyclic convex $n$-gon circumscribed in
a circle of radius $R$ with length of sides $a_1, a_2, a_3,...., a_n,\  L_n, P_n$ are respectively the perimeter and the pseudo-area. Then under hypothesis of Theorem 3.3, the following inequalities hold 
$$ 2 R \tan \frac {\pi }{n} \left(L_n -2n R \ sin \frac {\pi }{n} \right)\leq   {L_{{n}}}^{2}- \frac{4 P_{{n}}}{ \left( 1-\frac{2}{n} \right)^{\frac{n}{2}}}.$$
Equality holds if and only if  $\Pi_n $ is a regular polygon}.\\

This is a consequence of Proposition 4.2 and \cite{zd}, Corollary 6.\\

{\bf Corollary 4.4} \quad {\it Let $\Pi_n $ be a cyclic convex $n$-gon circumscribed in
a circle of radius $R$ with length of sides $a_1, a_2, a_3,...., a_n,\ A_n, L_n, P_n$ are respectively the area the perimeter and the pseudo-area. Then under hypothesis of Theorem 3.3, the following inequalities hold
$$ {L_{{n}}}^{2}- \frac{4 P_{{n}}}{ \left( 1-\frac{2}{n} \right)^{\frac{n}{2}}} \geq \frac{1}{R^2} \left[A_n - n R \left(sin \frac {2\pi }{n}\right)\right]^2.$$
 Equality holds if and only if  $\Pi_n $ is a regular polygon}.\\

This is a consequence of Proposition 4.2 and \cite{zd}, Theorem 9.\\

Of course, other inequalities can be deduced where $P_n$ appears by drawing inspiration from the results of \cite{zd}.\\

\section{Remarks and conclusions}

{\bf 1 - }\quad An interesting question that arises is whether it is possible to improve Theorems 3.1 and 3.3 by removing the assumption regarding the lengths of the sides.
I believe the answer is yes. To be convinced of this, one need only look at the multitude of examples provided by \cite{l} and \cite{ch2}. Indeed, let us consider, for instance the Macnab polygon, which is a cyclic 
alternate-sided $2n$-gon with $n$ sides of length $a$ and $n$ sides of length $b$. One proved \cite{ch2} the quotient $\varphi_{n,n}= \frac{A_{n,n}}{P_{n,n}}$ may be written
$$\varphi _{n,n}=\frac{[\cos \frac{\pi }{n}+\frac{2ab}{(a+b)^{2}}(1-\cos 
\frac{\pi }{n}]}{n\sin \frac{\pi }{n}[1-\frac{2}{n}+\frac{4ab}{%
n^{2}(a+b)^{2}}]^{\frac{n}{2}}}.$$ In addition, this quotient verify Conjecture (L)
$$\varphi_{n,n}^0 \leq \varphi_{n,n} \leq 1.$$
However, certain values ​​for the lengths $a,b$ may fail to satisfy the hypotheses of Theorems 3.1 and 3.3. Indeed, $$\ln  \left( \pi  \right) < n \frac{a}{L_{n,n}} = 2\,{\frac {na}{na+nb}}\quad \Rightarrow \quad\ln  \left( \pi  \right) < 2\,{\frac {a}{a+b}}.$$ To do this, we simply need to take \ $b < \frac{(2 -\log \pi)}{\log \pi} a.$\\

 {\bf 2 - }\quad In a 2-dimensional Euclidean space, a domain has an area if it is a Jordan-measurable set, and its area is equal to this measure.\\
In geometry and measurement theory, a standard area must respect three main axioms: always be positive, be zero for a point or a line, and be additive (the area of ​​two joined shapes is the sum of their areas). By relaxing or adapting one of these rules, we obtain different relevant definitions of a pseudo-area.\\
In short, for any planar polygon, one can define a measure called a "pseudo-area." This is a positive quantity that evaluates to zero for a point or a line, but it does not satisfy the additivity condition. In the last case, The area of ​​two combined shapes can be greater or smaller than the sum of the two (very useful in probability or decision theory). \\
If a polygon is deformed by "dis-articulating" it—while maintaining the same number of sides and the same side lengths— (that is to say, by deforming only one or more angles) then $P_n$ remains unchanged, but $A_n$ will decrease. Thus, the quotient $\phi_n = \frac{A_n}{P_n}$ will also decrease, and it can even approach \ $0.$\ Thus, for this non-cyclic $n$-gon its quotient is such that \ $\phi_n \leq \phi_n^0.$\\ 

The real question that arises is whether the inequality $\phi_n^0 \leq \phi_n$ for a $n$-gon $\Pi_n$ serves as an indicator for $\Pi_n$ to be cyclic. In other words, does the converse of Theorem 2.1 hold ?

$$\Pi_n \quad non-cyclic  \quad \Rightarrow \quad  \phi_n^0 = {\frac{1}{{{n \tan{\frac{\pi}{{n}}}}\ (1-{\frac{2}{{n}}}
)^{n/2}}}} \geq \phi_n \ ?$$

Recall the following result proved by \cite{o2}\\

{\bf Proposition 5.1 }\ {\it Let $K$ be a domain bounded by a rectifiable Jordan curve ${\cal C}$ with $L, A, R, r$\ its length,
area, circumradius and inradius. Then there exists a sequence of polygons $\Pi_n$ that
$$A n \rightarrow A,\quad L n \rightarrow L,\quad and\quad R_n   \rightarrow R  .$$
If in addition $K$ is convex, we may also arrange that\quad  $r_n \rightarrow r$  .}\\

 That means, intuitively, we may approximate ${\cal C}$ by finer polygonal curves.\\
What about the pseudo-area $P_n$ ? Suppose for a rectifiable Jordan curve ${\cal C}, \quad P_n \rightarrow P$\ exists. As for the area one has \ $A \leq \frac{L^2}{4 \pi}$\ then necessarily \ $P \leq \frac{L^2}{4 e}.$\ Indeed, we have seen before \cite{ch2} that for any $n$-gon (not necessary cyclic) 
$$\frac{4 P_n}{L_n^2 \left( 1-\frac{2}{n} \right)^{\frac{n}{2}}} \leq 1.$$
This means that when $n$ tends to $\infty$\ it yields
$$\lim_{\infty} \frac{4 P_n}{L_n^2 \left( 1-\frac{2}{n} \right)^{\frac{n}{2}}} = \frac{4 P e}{L^2} \leq 1.$$

\end{document}